\documentclass[reqno,draft,12pt]{amsart}

\usepackage[all]{xypic}
\usepackage{amsmath}
\usepackage{amscd}
\usepackage{amssymb}
\usepackage{enumerate}
\usepackage{amsfonts}
\usepackage{graphicx}
\usepackage[all]{xy}
\usepackage{mathrsfs}

\DeclareFontEncoding{OT2}{}{}

\numberwithin{equation}{section}

\newcommand{\rar}{\longrightarrow}
\newcommand{\hrar}{\hookrightarrow}
\newcommand{\trar}{\twoheadrightarrow}
\newcommand{\rarab}[1]{\overset{#1}{\longrightarrow}}
\newcommand{\hrarab}[1]{\overset{#1}{\hookrightarrow}}

\newcommand{\al}{\alpha}

\newcommand{\ga}{\gamma}
\newcommand{\Ga}{\Gamma}
\newcommand{\de}{\delta}
\newcommand{\De}{\Delta}
\newcommand{\la}{\lambda}
\newcommand{\La}{\Lambda}
\newcommand{\ze}{\zeta}

\newcommand{\eps}{\epsilon}
\newcommand{\sg}{\sigma}

\newcommand{\om}{\omega}
\newcommand{\Om}{\Omega}
\newcommand{\vp}{\varpi}

\newcommand{\ob}{\overline{\omega}}

\newcommand{\bA}{{\mathbb A}}

\newcommand{\bC}{{\mathbb C}}

\newcommand{\bG}{{\mathbb G}}

\newcommand{\bN}{{\mathbb N}}

\newcommand{\bQ}{{\mathbb Q}}
\newcommand{\bR}{{\mathbb R}}

\newcommand{\bZ}{{\mathbb Z}}

\newcommand{\cN}{{\mathcal N}}

\newcommand{\fa}{{\mathfrak a}}

\newcommand{\fg}{{\mathfrak g}}

\newcommand{\Rt}{\widetilde{R}}
\newcommand{\Vt}{\widetilde{V}}
\newcommand{\alt}{\widetilde{\al}}

\newcommand{\mb}{\overline{m}}

\newcommand{\Hom}{\operatorname{Hom}}

\newcommand{\Lie}{\operatorname{Lie}}
\newcommand{\Spec}{\operatorname{Spec}}

\newcommand{\pr}{\mathrm{pr}}

\newcommand{\Gal}{\operatorname{Gal}}

\newcommand{\val}{\operatorname{val}}
\newcommand{\tens}{\otimes}

\newcommand{\Res}{\operatorname{Res}}

\newcommand{\sbr}{\smallbreak}
\newcommand{\mbr}{\medbreak}

\newcommand{\indlim}[1]{\lim\limits_{\underset{N}{\longrightarrow}}}
\newtheorem{thm}{Theorem}[section]
\newtheorem{cor}[thm]{Corollary}
\newtheorem{lem}[thm]{Lemma}

\newtheorem{prop}[thm]{Proposition}

\theoremstyle{remark}
\newtheorem{rem}[thm]{Remark}

\newcommand{\Fb}{\overline{F}}

\newcommand{\bbe}{\begin{equation}}
\newcommand{\ee}{\end{equation}}

\newcommand{\pair}[1]{\langle #1\rangle}

\newcommand{\ssect}{\subsection}

\newcommand{\cst}{\bC^\times}

\newcommand{\Ker}{\operatorname{Ker}}

\newcommand{\Exp}{\operatorname{Exp}}

\title[Centralizers of generic elements of Newton strata]{Centralizers
of generic elements of Newton strata in the adjoint quotients of reductive
groups}
\author[Mitya Boyarchenko and Maria Sabitova]{Mitya
Boyarchenko\address{Mitya Boyarchenko: University of Chicago, \hfill\newline
Department of Mathematics, Chicago, IL 60637} and Maria Sabitova\address{Maria
Sabitova: University of Illinois at Urbana-Champaign, \hfill\newline Department
of Mathematics, Urbana, IL 61801}}

\thanks{The research of M.B. was partially supported by NSF grant
DMS-0401164.  \hfill\newline {\em Email addresses:}\ \ {\tt
mitya@math.uchicago.edu} (M.B.), \ \ {\tt sabitova@math.upenn.edu} (M.S.)}

\begin{document}

\begin{abstract}
We study the Newton stratification of the adjoint quotient of a connected split
reductive group $G$ with simply connected derived group over the field $F=\bC((\eps))$ of formal Laurent series. Our main result describes the
centralizer of a regular semisimple element in $G(F)$ whose image in the
adjoint quotient lies in a certain generic subset of a given Newton stratum.
Other noteworthy results include analogues of some results of Springer on
regular elements of finite reflection groups, as well as a geometric
construction of a well known homomorphism $\psi:P(R^\vee)/Q(R^\vee)\to W(R)$ defined
for every reduced and irreducible root system $R$.
\end{abstract}

\maketitle

\setcounter{tocdepth}{1}

\tableofcontents


\section*{Introduction}

If $F=\bC((\eps))$ is the field of formal Laurent series and $G=GL_n(\Fb)$, it
is well known that a semisimple element of $G$ is conjugate to an element lying
in $GL_n(F)$ if and only if the coefficients of its characteristic polynomial
lie in $F$, and, conversely, every monic polynomial of degree $n$ over $F$
whose constant term in nonzero is the characteristic polynomial of a semisimple
element of $G$ lying in $GL_n(F)$. Hence the $F$-points of the adjoint quotient
$\bA$ of $G$ are identified with the set  of such
polynomials: $\bA(F)\cong F^{n-1}\times F^\times$. We can partition $\bA(F)$ into a disjoint union of subsets
$\bA(F)_{\nu}$ consisting of polynomials of $\bA(F)$ having the same Newton
polygon $\nu$, which gives the {\em Newton stratification} of $\bA(F)$.

\mbr

In his paper \cite{k} R.~Kottwitz generalizes the notion of Newton
stratification to a general connected split reductive group $G$ over $F$ whose
derived group is simply connected. One of the main goals of this paper is to
answer a question (also formulated by Kottwitz) which arises naturally in this
context. Namely, we prove that for every Newton stratum
$\bA(F)_\nu\subset\bA(F)$ and for every regular semisimple element $\ga\in
G(F)$ whose image in $\bA(F)$ lies within a certain open dense subset of
$\bA(F)_\nu$, the centralizer $G_\ga$ is isomorphic to the twist of a split
maximal torus $A\subset G$ by an element $w=w(\nu)$ of the Weyl group of $G$
with respect to $A$. The element $w$ is independent of the choice of $\ga$ and admits a
rather simple description in terms of $\nu$ (see Theorem \ref{t:1}). The
definition of the Newton stratification appears in Section \ref{s:notation}
together with some other preliminaries, while the statement and the proof of
our main result occupy Sections \ref{s:main} and \ref{s:systems}.

\mbr

In the course of the proof we encountered several questions about semisimple
groups over $\bC$ that, to the best of our knowledge, have not been studied
before. In particular, if $G$ is a simple, connected and simply connected
algebraic group over $\bC$ with maximal torus $A\subset G$, we establish in
Section \ref{s:Springer} natural analogues of some of the results of
T.~Springer \cite{s} in the nonlinear context of the Weyl group $W=W(G,A)$
acting on $A$. More precisely, we prove that, given an element $x$ of the center $Z(G)$ of $G$,
there exist $g\in W$ and a regular element $u\in A$ satisfying $g(u)=xu$; in
addition, this condition determines the element $g$ uniquely up to conjugation,
and the eigenvalues of the action of $g$ on the Lie algebra $\Lie(A)$ can be
computed easily in terms of $x$. Aside from the applications appearing in this
paper, the usefulness of these results is demonstrated by Remark
\ref{r:kottwitz}.

\mbr

Section \ref{s:geometric} is devoted to the study of a certain important
homomorphism $\psi:Z(G)\to W$ (which already enters the formulation of our
first main result). This homomorphism is well defined up to conjugation by
elements of $W$. So far the only known (to us) description of this homomorphism
was purely algebraic. Our last main result is Theorem \ref{t:monodromy}, which
provides a geometric (or, if the reader wishes, topological) description of
$\psi$ as the monodromy of a certain $W$-torsor over a manifold whose
fundamental group is identified with $Z(G)$, constructed by restricting the
Springer covering $\widetilde{G}^{rs}\to G^{rs}$ to the set of regular elements
of a maximal compact subgroup of $G$.

\subsection*{Acknowledgements} We would like to express our deepest
gratitude to Robert Kottwitz for introducing us to this subject and suggesting
the statement of one of our main results, as well as carefully reading our article. We would also like to thank him for
sharing the preprint \cite{k} before it was made available to the general
audience, for his time, for numerous stimulating discussions, and in particular
for pointing out a way of greatly simplifying the proof of the first part of
Theorem \ref{t:monodromy}. We are grateful to Misha Finkelberg for useful and
inspiring conversations, for referring us to the paper \cite{ss}, and for
explaining to us Zakharevich's construction that appears in
\S\ref{ss:monodromy}. We are indebted to Jim Humphreys, David Kazhdan, George
Lusztig, Jean-Pierre Serre, and Tonny Springer for very helpful e-mail
correspondence. Part of the research was conducted when the second author was a
postdoc at MSRI in the Spring of 2006, and she would like to thank the
Institute for its hospitality.


\section{General facts and notation}\label{s:notation}

Since the main reference for the first part of this paper is \cite{k}, we will
keep the same notation as in that article. In this section we review this
notation and recall some results that will be used in the remainder of the
paper. The reader familiar with \cite{k} may wish to skip directly to Section
\ref{s:main} and refer back to this one as the need arises.

\subsection{The setup}\label{ss:setup} Let $G$ be a connected reductive group
over $\bC$ of rank $n$. We assume that the derived group $G_{sc}=[G,G]$ of $G$
is simply connected, and let $l$ be its rank (so that $l\leq n$). The principal
example to keep in mind is $G=GL_n$. We fix a maximal torus $A\subset G$ and a
Borel subgroup $B\subset G$ containing $A$, and write $A_{sc}=A\cap G_{sc}$,
which is a maximal torus of $G_{sc}$. We let $R\subset X^*(A)$ denote the root
system of $G$ with respect to $A$, and $\De=\{\al_i\}_{i\in I}\subset R$ the
basis of simple roots determined by $B$, where $I=\{1,2,\dotsc,l\}$. Further,
we let $Q(R)$ and $Q(R^\vee)$ denote the root lattice and the coroot lattice of
$R$, and we let $P(R)$ and $P(R^\vee)$ denote the dual lattices to $Q(R^\vee)$
and $Q(R)$, respectively. It is known that $P(R)$ can be canonically identified
with $X^*(A_{sc})$, and $Q(R^\vee)$ can be canonically identified with
$X_*(A_{sc})$. We will view these identifications as equalities. Finally, let
$\{\vp_i\}_{i\in I}$ denote the basis of fundamental weights of
$P(R)=X^*(A_{sc})$ which is dual to the basis of $Q(R^\vee)$ consisting of the
simple coroots $\{\al_i^\vee\}_{i\in I}$, and let $\om_1,\ldots,\om_n\in
X^*(A)$ be a $\bZ$-basis chosen as in \cite{k} (\S 1.3, p.~3). Namely,
$\om_1,\ldots,\om_l$ are preimages of the fundamental weights
$\varpi_1,\ldots,\varpi_l\in X^*(A_{sc})$ under the natural map
$\Res:X^*(A)\twoheadrightarrow X^*(A_{sc})$ and $\om_{l+1},\ldots,\om_n$ is a
$\bZ$-basis of the kernel of $\Res$.

\subsection{The adjoint quotient of $G$}\label{ss:adj-quotient}
Let $W=W(G,A)=N_G(A)/A$ be the Weyl group of $G$ with respect to $A$. For each
$\la\in X^*(A)$ we write $e^{\la}$ for $\la$ viewed as an element of the group
algebra $\bC[X^*(A)]$ of $X^*(A)$ over $\bC$; note that
$\bC[X^*(A)]\cong\bC[A]$, the coordinate ring of $A$. Put
$$
c_i:=\sum_{\la\in W\om_i}e^{\la}\in \bC[X^*(A)]^W,\quad i\in\{1,\ldots,n \},
$$
where $W\om_i$ denotes the $W$-orbit of $\om_i$ in $X^*(A)$. Let
$\bA:=\bG_a^l\times\bG_m^{n-l}$. It is well-known that the map $c:A\rar\bA$ defined
by $$a\mapsto (c_1(a),\ldots,c_n(a))$$ induces an isomorphism of algebraic
varieties between $A/W$ and $\bA$ (see \cite{k}, \S1.5, p. 5, for example). We will often
identify $A/W$ with $\bA$ via $c$. Note that, in particular, we obtain an
identification $A(\overline{F})/W\cong\bA(\overline{F})$ for any algebraically
closed field $\overline{F}$ containing $\bC$.

\mbr

The variety $\bA$ is called the {\em adjoint quotient} of $G$, because if we
let $G$ act on itself by conjugation, Chevalley's restriction theorem implies
that the restriction map $\bC[G]\to\bC[A]$ induces an isomorphism of algebras
$\bC[G]^G\cong\bC[A]^W$, and hence we obtain an isomorphism of algebraic
varieties $A/W\cong G/\!/(\operatorname{Ad} G)$.

\subsection{The base field $F$}\label{ss:base-field}
Let $F=\bC((\eps))$ be the field of formal Laurent power series
over $\bC$ in an indeterminate $\eps$.  For each $n\in\bN$ we let $\eps^{1/n}\in\Fb^{\times}$
be a fixed $n$-th root of $\eps$
such that
$$
\left(\eps^{\frac{1}{mn}}\right)^{m}=\eps^{1/n},\quad
\forall m,n\in\bN.
$$
It is known that $\Fb=\bigcup_{n\in\bN}\bC((\eps^{1/n}))$. Thus, the element
$\sg\in\Ga=\Gal(\Fb/F)$ given by
$$
\sg(\eps^{1/n})=\exp\left(\frac{2\pi
\sqrt{-1}}{n}\right)\eps^{1/n},\quad\forall n\in\bN,
$$
is a topological generator of $\Ga$, i.e., it determines an
isomorphism $\widehat{\bZ}\rarab{\simeq}\Ga$. We also fix the valuation $\val$ on $\Fb$ such that
$$
\val(\eps^{1/n})=-1/n,\quad\forall n\in\bN, \quad\text{ and }\,\,\val(0)=-\infty.
$$

\mbr

We may view $G$ and $A$ as algebraic groups over $F$ by extending scalars. This
should cause no confusion, since in the remainder of this section, as well as
in Sections \ref{s:main} and \ref{s:systems} we always think of $G$ and $A$ as
algebraic groups over $F$, whereas in Sections \ref{s:Springer} and
\ref{s:geometric} we will only work with algebraic groups over $\bC$. We will
consider $\bA(F)=F^l\times (F^\times)^{n-l}$ as a topological space with
respect to the topology induced from $F^n$, where $F$ is identified with a
product of countably many copies of $\bC$ and the closed sets in $F^n$ are
defined by polynomial equations over $\bC$ in finitely many coordinates.

\subsection{The map $f:G^{rs}(F)\rar\bA(F)$}\label{subsec:f}
Let $G^{rs}(F)$ denote the set of regular semisimple elements of $G$ defined
over $F$. We will need the map $f:G^{rs}(F)\rar\bA(F)$ defined in the following
way. Let $\ga\in G^{rs}(F)$ be a regular semisimple element of $G$, i.e.,
$\ga\in G(F)$ and the centralizer $G_{\ga}$ of $\ga$ in $G$ is a maximal torus
(defined over $F$). There is $g\in G(\Fb)$ which conjugates $\ga$ to an element
$\ga_0\in A(\Fb)$. Note that $\ga_0$ is only determined up to the $W$-action.
Namely, if $g_1\ga g_1^{-1}\in A(\Fb)$ and $g_2\ga g_2^{-1}\in A(\Fb)$, then
the element $g_1g_2^{-1}\in G(\Fb)$ conjugates $g_2\ga g_2^{-1}$ into $g_1\ga
g_1^{-1}$. Since $\ga$ is regular and semisimple,  we have
$$
A(\Fb)=G_{g_1\ga g_1^{-1}}(\Fb)=G_{g_2\ga g_2^{-1}}(\Fb),
$$
hence $g_1g_2^{-1}\in N_G(A)(\Fb)$, where $N_G(A)$ denotes the normalizer of
$A$ in $G$. Thus, $g_1\ga g_1^{-1}$ and $g_2\ga g_2^{-1}$ are in the same
$W$-orbit. So we get a well-defined map $f:G^{rs}(F)\rar\bA(\Fb)$ given by
$$
f(\ga)=c(\ga_0).
$$
We claim that, in fact, $f(\ga)\in\bA (F)$. To check the claim we need to show
that the $W$-orbit of $\ga_0$ in $A(\Fb)$ is stable under the action of $\Ga$
on $A(\Fb)$. Let $\la\in\Ga$, then
$$
\la(\ga_0)=\la(g\ga g^{-1})=\la(g)\ga\la(g)^{-1},
$$
because $\ga\in G(F)$. This computation shows that $\la(\ga_0)$ is conjugate to
$\ga_0$ by an element of $G(\Fb)$ and hence, by the argument above,
$\la(\ga_0)$ and $\ga_0$ lie in the same $W$-orbit, hence the image of $f$ lies
in $\bA(F)$.

\subsection{The homomorphism $\psi_G:\La_G\rar W$}\label{ss:psi-G}
Let $\La_G=X_*(A)/X_*(A_{sc})$, and let $A_G\subseteq A$ be the identity
component of the center of $G$. We have the following natural maps:
\bbe\label{eq:l} \La_G\trar\La_G/X_*(A_G)\hrar P(R^{\vee})/Q(R^{\vee}), \ee
where $X_*(A_G)$ is identified with its image under the embedding
$X_*(A_G)\hrar X_*(A)\trar \La_G$ and $R^{\vee}$ denotes the coroot system of
$G$. Let $\psi:P(R^{\vee})/Q(R^{\vee})\rar W$ be the map defined as in \cite{b}
(Chapter VI, \S2, no.~3). (See also \S\ref{ss:psi} below.) We denote by
$\psi_G$ the composition
$$
\La_G\rar P(R^{\vee})/Q(R^{\vee})\rarab{\psi} W,
$$
where the first map is given by \eqref{eq:l}.

\begin{rem} We note that $\psi_G$ is the map $ps$ introduced
in \cite{k} (\S 1.9, p. 9).
\end{rem}

\subsection{Newton stratification of $\bA(F)$} Let
$\fa=X_*(A)\otimes_{\bZ}\bR$ and let $P=MN$ be a parabolic subgroup of $G$
containing $B$ (i.e., $P$ is a {\em standard} parabolic subgroup) with  Levi
subgroup $M$ and  unipotent radical $N$. Let $W_M$ be the Weyl group of $M$
identified with a subgroup of $W$. Put $\fa_M:=X_*(A_M)\otimes_{\bZ}\bR$, where
as above $A_M\subseteq A$ denotes the identity component of the center of $M$.
Then $\fa_M$ can be identified with the set of fixed points of $\fa$ under the
action of $W_M$, and we have a surjection
$$
p_M:\fa\trar \fa_M,\quad x\mapsto\frac{1}{W_M}\cdot\sum_{w\in W_M}w(x).
$$
Since $G_{sc}$ is simply connected, the derived group of $M$ is also simply
connected, which implies that $\La_M$ can be identified with the image of
$X_*(A)$ under $p_M$. In what follows we write $\La_M\hrar\fa_M$ to indicate
that $\La_M$ is considered as  $p_M(X_*(A))\subseteq\fa_M$.

\mbr

Let $\cN_G$ be the subset of $\fa$ defined as in \cite{k} (\S 1.3, p. 4). For the
sake of completeness we repeat the facts about $\cN_G$ which will be used
in this paper. First,
$$
\cN_G=\coprod_{P}\La_P^+ \subseteq\fa,
$$
where the union is taken over all standard parabolic subgroups of $G$ and for
each such $P$ (in the above notation) $\La_P^+$ is a subset of
$\La_M\hrar\fa_M$. Moreover, $\La_G^+=\La_G$, where $\La_G\hrar\fa_G$. Thus,
for each $\nu\in\cN_G$ there is a unique parabolic subgroup $P=MN$ such that
$\nu\in\La_P^+$, hence the element $\psi_M(\nu)\in W_M$ is defined and we put
\bbe\label{eq:we} w(\nu):=\psi_M(\nu)\in W.\ee

Also, for each $\nu\in\cN_G$ there is a certain non-empty irreducible subset
$\bA(F)_{\nu}$ of $\bA(F)$ such that
$$
\bA(F)=\coprod_{\nu\in\cN_G}\bA(F)_{\nu}
$$
(\cite{k}, Thm. 1.5.2, p. 6).

\mbr

In what follows we will need one result about the sets $\bA(F)_{\nu}$, which we
state as a lemma for future reference.

\begin{lem}\label{l:kot}
Let $\nu\in\cN_G$. Then $c\in\bA(\Fb)$ belongs to $\bA(F)_{\nu}$ if and only if
there is $a\in A(\Fb)$ such that
\begin{eqnarray}
&& c=c(a),\label{eq:surj1}\\
&&\sg(a)=g(a)\quad\text{for some }g\in W,\text{ and}\label{eq:surj2}\\
&& \pair{\la,\nu}=\val\la(a)\quad\text{for any }\la\in X^*(A).\label{eq:surj3}
\end{eqnarray}
\end{lem}

\begin{proof}
See Theorem 1.5.2(5) of \cite{k}.
\end{proof}

\mbr

We consider $\bA(F)_{\nu}$ endowed with the topology induced from the topology
on $\bA(F)$ described in \S\ref{ss:base-field}. By a {\em generic subset} of
$\bA(F)_{\nu}$ we mean a subset containing an open dense subset of
$\bA(F)_{\nu}$. Also, we say that a certain property {\em holds for a generic
element of $\bA(F)_{\nu}$} if it holds for all elements of some generic subset
of $\bA(F)_{\nu}$.

\subsection{Twists of a torus}\label{ss:tor} Let $w\in W$. In what follows we
will denote by $A^w$ the torus of $G$ obtained from $A$ by twisting the Galois
structure on $A$ by $w$. Then $X^*(A^w)=X^*(A)$ as abelian groups, and if
$\la\in X^*(A^w)$, then the Galois action on $\la$ is given by
$$
\sg\la(a)=\la(w^{-1}(a)), \quad a\in A.
$$
We say that a torus $T$ of $G$ defined over $F$ is {\em obtained from $A$ by
twisting by $w$} if and only if $X^*(T)\cong X^*(A^w)$ as $\Ga$-modules (with
the natural action of $\Ga$ on $X^*(T)$).


\section{Statement of the main theorem}\label{s:main}

We keep all the notation introduced above. In this section we view $G$ and $A$
as algebraic groups over $F$. Our goal is to state the main result mentioned in
the introduction (Theorem \ref{t:1}), give a more precise reformulation of this
result (Theorem \ref{t:2}), and reduce its proof to a special case (Theorem
\ref{t:3}).

\subsection{Main result} The next two sections, as well as \S\ref{ss:Conclusion},
are devoted to the proof of

\begin{thm}\label{t:1}
Let $\nu\in\cN_G.$ Then for generic $c\in\bA(F)_{\nu}$  and a regular
semisimple element $\ga\in G(F)$ such that $f(\ga)=c$ we have
$$
G_{\ga}\cong A^{w(\nu)},
$$
where $w(\nu)$ is given by \eqref{eq:we}.
\end{thm}

Let us slightly reformulate Theorem \ref{t:1}. As we explained in
\S\ref{subsec:f}, there exists $g\in G(\Fb)$ such that
\begin{eqnarray*}
g\ga g^{-1}&=&\ga_0\in A(\Fb)\quad\text{and}\\
g\sg(g)^{-1}&\in &N_G(A)(\Fb).
\end{eqnarray*}
Let $h=h(\ga)$ denote the element of $W$ corresponding to $g\sg(g)^{-1}$. Note
that $h$ is defined up to conjugation in $W$.

\begin{lem}
We have $G_{\ga}\cong A^h.$
\end{lem}
\begin{proof} Since $gG_{\ga}g^{-1}=A$, conjugation by $g$ induces
the isomorphism $\phi:X^*(A)\rar X^*(G_{\ga})$. By \S\ref{ss:tor} it is enough
to show that $\phi$ commutes with the action of $\sg$. Equivalently,
\bbe\label{eq:gx} (\sg\la)(gxg^{-1})=\sg\big(\phi(\la)\big)(x),\quad\text{for
any }\la\in X^*(A),\,x\in G_{\ga}.\ee Since $A$ splits, the right side of
\eqref{eq:gx} is equal to $\la\big(\sg(g)x\sg(g)^{-1}\big)$, and by the
definition of the Galois action on $X^*(A^h)$, the left hand side of \eqref{eq:gx} is equal to
$$\la\big(h^{-1}(gxg^{-1})\big)=\la\big(\sg(g)x\sg(g)^{-1}\big).$$
\end{proof}

According to this lemma to prove Theorem \ref{t:1} it is enough to prove the
following

\begin{thm}\label{t:2'}
Let $\nu\in\cN_G$. There exists a generic subset $X_\nu\subset\bA(F)_\nu$ such
that given $a\in A(\Fb)$ for which $G_a$ is a torus and $c(a)\in X_{\nu}$, the
element $h\in W$ given by \bbe \sg(a)=h(a)\ee is conjugate to $w(\nu)$.
\end{thm}

\subsection{Construction of $X_\nu$}
Recall that $R\subset X^*(A)$ denotes the root system of $G$, let $R^+\subset
R$ be the set of positive roots with respect to $B$, and let
$$
\Om:=\prod_{\al\in R}\big(1-e^{\al}\big)\in \bC[X^*(A)].
$$
For $\nu\in\cN_G$ denote
$$
m_{\nu}:=\sum_{\al\in R}\max\left\{0,\pair{\al,\nu}\right\}.
$$
Observe that $m_{\nu}=\pair{2\rho,\nu}$, where $$\rho=\frac{1}{2}\cdot\sum_
{\al\in R^{+}}\al\in X^*(A).$$ (This follows from the facts mentioned in the
first section of the proof of Lemma \ref{l:en} below.) Also, if
$\nu\in\La_G\hrar\fa_G$ then $m_{\nu}=0$. Note that for any $a\in A$ such that
$c(a)\in\bA(F)_{\nu}$, we have $\val\Om(a)\leq m_{\nu}$. Indeed, by Lemma
\ref{l:kot} there exists $b\in A$ such that $c(a)=c(b)$ and
$\val\la(b)=\pair{\la,\nu}$ for any $\la\in X^*(A)$. Hence $a=g(b)$ for some
$g\in W$, and consequently
$$
\val\Om(a)=\val\Om(b)=\sum_{\al\in R}\val(1-\al(b))\leq\sum_{\al\in
R}\max\left\{0,\pair{\al,\nu}\right\}=m_{\nu}.
$$
(Recall that with our conventions, $\val$ is the negative of the usual valuation, so that $$\val(a+b)\leq\max\{\val a,\val b\},\quad\text{for any }a,b\in\Fb.)$$
On the other hand, $\Om\in\bC[X^*(A)]$ belongs to the subring $\bC[X^*(A)]^W$
of $W$-invariant elements in $\bC[X^*(A)]$. Since $G$ comes from an algebraic
group defined over $\bC$,
$$
\bC[X^*(A)]^W=\bC[c_1,\ldots,c_n,c_{l+1}^{-1},\ldots,c_{n}^{-1}].
$$
Together with the note above, this implies that
$$
X_{\nu}:=\big\{c\in\bA(F)_{\nu}\,\big\vert\,\exists\, a\in A\text{ s.t. }
c=c(a)\text{ and }\val\Om(a)=m_{\nu}\big\}
$$
is an open subset of $\bA(F)_{\nu}$.

\mbr

Since $\bA(F)_{\nu}$ is irreducible, Theorem \ref{t:2'} is a consequence of the
following result:
\begin{thm}\label{t:2}
For each $\nu\in\cN_G$ the set $X_{\nu}$ is non-empty and for any $a\in A$ such
that $c(a)\in X_{\nu}$, the element $h\in W$ given by \bbe \sg(a)=h(a) \ee is
conjugate to $w(\nu)$.
\end{thm}

\begin{rem}\label{r:chert}
Observe that the condition $c(a)\in X_{\nu}$ implies automatically that $G_a$
is torus, or equivalently, $G_a=A$. Indeed, write $a=a_1a_2$, where $a_1\in
A_{sc}$ and $a_2\in A_G$. Then $G_a=(G_{sc})_{a_1}\cdot A_G$. Since $a_1\in
A_{sc}$, it is in particular a semisimple element of $G_{sc}$, and hence the
centralizer $(G_{sc})_{a_1}$ of $a_1$ in $G_{sc}$ is {\em connected} by a deep
result of Springer and Steinberg (\cite{ss}, Theorem II.3.9), which is false
without the assumption that $G_{sc}$ is simply connected. Thus, $G_a$ is
connected, as a product of two connected algebraic groups. Now, assume that
$G_a$ does not coincide with $A$. Then the Lie algebra $\Lie(A)$ of $A$ is a
proper subspace of the Lie algebra $\Lie(G_a)$ of $G_a$. Moreover, $\Lie(G_a)$
is invariant under the adjoint action of $A$, since $A\subset G_a$. This
implies that $\Lie(G_a)$ contains a root subspace $\fg_{\al}$ for some $\al\in R$. Thus
$\al(a)=1$, and we get a contradiction with $c(a)\in X_{\nu}$.
\end{rem}

\subsection{Reduction from $\cN_G$ to $\La_G$}
To end this section, we observe that in Theorem \ref{t:2} it is enough to
assume that $\nu\in\La_G\hrar\fa_G$, which means that we only have to prove

\begin{thm}\label{t:3}
For each $\nu\in\La_G\hrar\fa_G$ the set $X_{\nu}$ is non-empty and for any
$a\in A$ such that $c(a)\in X_{\nu}$, the element $h\in W$ given by \bbe
\sg(a)=h(a) \ee is conjugate to $w(\nu)=\psi_G(\nu)$.
\end{thm}

Indeed, we have the following result.

\begin{lem}\label{l:en}
Theorem $\ref{t:3}$ implies Theorem $\ref{t:2}.$
\end{lem}

\begin{proof}
Let $\nu\in\La_P^+$ for some standard parabolic subgroup $P$ of $G$ with the
Levi subgroup $M$ containing $A$, so that $\La_P^{+}\subseteq\La_M\hrar\fa_M$.
Denote by $\De_M\subseteq\De$ the set of simple roots of $M$. Also, let
$R_M\subseteq R$ denote the root system of $M$ and let $$\Om_M:=\prod_{\al\in
R_M}(1-e^{\al})\in\bC[X^*(A)].$$ Recall (\cite{k}, p. 4) that $\La_P^+$ is contained in
the set of elements
 $x\in\fa_M$ such that $\pair{\al,x}>0$ for any $\al\in\De\backslash\De_M$.
Since
$$\fa_M=\big\{x\in\fa\,\big\vert\,\pair{\al,x}=0,\,\forall\al\in \De_M\big\},
$$
we conclude that $\pair{\al,\nu}=0$ for any $\al\in R_M$ and $\pair{\al,\nu}\ne
0$ for any $\al\in R\backslash R_M$. Thus, if $a\in A$ satisfies
\eqref{eq:surj3} and $c(a)\in\bA(F)_{\nu}$, then
\bbe\label{eq:op}\val\Om(a)=m_{\nu}\Longleftrightarrow\val\Om_M(a)=0.\ee

\mbr

Let $\bA_M$ and $c_M:A/W_M\rarab{\sim}\bA_M$ denote respectively the set $\bA$
and the map $c$ corresponding to $M$. Consider the map from $\bA_M$ to $\bA$
defined by the composition of $c_M^{-1}$ with the natural map $A/W_M\rar A/W$
(induced by the embedding $W_M\hookrightarrow W$) followed by $c$. It is easy
to check that this map induces a map $\pi_M:\bA_M(F)_{\nu}\rar\bA(F)_{\nu}$. We
are now ready to prove the lemma. Let us show first that $X_{\nu}$ is
non-empty. Let $X_{\nu}(M)$ denote the set $X_{\nu}$ defined for $M$. By
Theorem \ref{t:3} and Lemma \ref{l:kot}, there exists $a\in A$ satisfying
\eqref{eq:surj3} such that $c_M(a)\in\bA_M(F)_{\nu}$ and $\val\Om_M(a)=0$. Then
$c(a)=\pi_M(c_M(a))\in\bA(F)_{\nu}$ and $\val\Om(a)=m_{\nu}$ by \eqref{eq:op}.
Thus, $c(a)\in X_{\nu}$.

\mbr

Now we will show that Theorem \ref{t:2} holds for $X_{\nu}$. Let $a\in A$
satisfy $c(a)\in X_{\nu}$ and $\sg(a)=h(a)$. Then by Lemma \ref{l:kot} there
exists $b\in A$ satisfying  \eqref{eq:surj3} such that $c(a)=c(b)$, hence
$b=s(a)$ for some $s\in W$. Thus, $\sg(b)=g(b)$ for some $g\in W$. Since $a$ (and hence $b$) is regular semisimple, as it follows
from Remark \ref{r:chert}, we see that $g$ is conjugate to $h$. Hence it is enough to
show that $g$ is conjugate to $w(\nu)=\psi_M(\nu)$. For any $\la\in X^*(A)$ we
have
$$
\pair{\la,\nu}=\val\la(b)=\val\sg(\la(b))=\val g^{-1}\la(b)=\pair{\la,g\nu},
$$
which implies $\nu=g\nu$. Thus, $g\in W_M$ by standard facts about reflection
groups, cf. \cite{b}. By Lemma \ref{l:kot} this gives $c_M(b)\in\bA_M(F)_{\nu}$
and also $\val\Om_M(b)=0$ by \eqref{eq:op}. Hence, $c_M(b)\in X_{\nu}(M)$ and
$g$ is conjugate to $\psi_M(\nu)$ by Theorem \ref{t:3}.
\end{proof}


\section{A result about root systems}\label{s:systems}

The goal of this section is to show that Theorem \ref{t:3} follows from a
certain statement about root systems (Proposition \ref{p:sys} below). This
provides a link between our main result (which is specific to algebraic groups
over $F$) and the theory of semisimple algebraic groups over $\bC$. Note that
the torus $T$ appearing in \S\ref{ss:generalities} corresponds to the torus
$A_{sc}$ defined in \S\ref{ss:setup}; otherwise our notation remains the same
as before.

\subsection{Generalities}\label{ss:generalities}
Let $R$ be a (reduced) root system in a real vector space $V$ and let $T$ be
the complex torus with character lattice $P(R)$, i.e.,
$$
T=\Spec\bC[P(R)].
$$
For each $\mu\in P(R^{\vee})$ we let $x_{\mu}$ denote the composition
$$
P(R)\rarab{\pair{\cdot,\mu}}\bQ/\bZ\hrarab{\exp(-2\pi
\sqrt{-1}\cdot)}\bC^{\times}.
$$
Thus $x_{\mu}$ is a group homomorphism from $P(R)$ to $\bC^{\times}$, and we
identify it with the (complex) point of $T$ it defines.

\mbr

Let $W$ denote the Weyl group of $R$, acting on $P(R)$ and hence on $T$ in  the
usual way, and let
$$
\psi:P(R^{\vee})/Q(R^{\vee})\rar W
$$
denote the homomorphism mentioned in \S\ref{ss:psi-G} above (see also \S\ref{ss:psi}).

\mbr

Recall that an element $u\in T$ is said to be {\em regular} if the stabilizer
of $u$ in $W$ is trivial, or, equivalently, if
$$
\al(u)\ne 1,\quad\forall\al\in R.
$$
(The equivalence of the two properties is not an obvious statement, and follows
from the same argument as in Remark \ref{r:chert}.)

\subsection{Reduction of Theorem \ref{t:3} to Proposition \ref{p:sys}}\label{ss:reduction}
In this section we assume

\begin{prop}\label{p:sys}
Let $\mu\in P(R^{\vee})$ and let $h\in W$. There exists a regular element
$u\in T$ with $h(u)=x_{\mu} u$ if and only if $h$ is conjugate to
$\psi(\bar{\mu})$.
\end{prop}

The proof is given in \S\ref{ss:Conclusion}.

\begin{lem}\label{l:enough}
Proposition $\ref{p:sys}$ implies Theorem $\ref{t:3}.$
\end{lem}

\begin{proof}
In the notation of Theorem \ref{t:3}, let $R$ be the root system of $G$ with
respect to $A$, and note that the torus $T=\Spec\bC[P(R)]$ is identified with the maximal
torus $A_{sc}$ of $G_{sc}$. Fix $\nu\in \La_G$. We first show that $X_{\nu}$ is
non-empty. Put
$$
m_i:=\pair{\om_i,\nu}\in\bQ,\quad i\in\{1,\ldots,n\},
$$
where as usual we consider $\nu$ as an element of $\fa_G$ under the embedding
$\La_G\hookrightarrow\fa_G$. Let $\mu\in P(R^{\vee})$ be such that
$\nu\mapsto\bar{\mu}$ under the map
$$
\La_G\twoheadrightarrow\La_G/X_*(A_G)\hrar P(R^{\vee})/Q(R^{\vee}).
$$
Then \bbe\label{eq:equiv}
\pair{\om_i,\nu}\equiv-\pair{\varpi_i,\mu}\mod\,\bZ\quad(1\leq i\leq l) \ee
(\cite{k}, p. 14). Furthermore, by definition $\psi_G(\nu)=\psi(\bar{\mu})$. By
Proposition \ref{p:sys} applied to $R$ and $\mu$ there exists $g\in W$ and a
regular element $u\in T$ such that $g(u)=x_{\mu} u$. Note that in view of
\eqref{eq:equiv} we have $g(u)=x_{\mu} u$ if and only if \bbe\label{eq:starr}
g^{-1}\varpi_i(u)=\exp\left(2\pi\sqrt{-1}m_i\right)\cdot\varpi_i(u),\quad
\forall i\in\{1,\ldots,l\}. \ee Let $a\in A$ be given by
\begin{eqnarray*}
\om_i(a)&=&\eps^{m_i}\cdot\varpi_i(u),\quad i\in\{1,\ldots,l\},\\
\om_i(a)&=&\eps^{m_i},\,\,\,\,\,\,\,\quad\quad\quad i\in\{l+1,\ldots,n\}.
\end{eqnarray*}
We claim that $c(a)\in X_{\nu}$. Recall that $c\in\bA(\Fb)$ belongs to
$X_{\nu}$ if and only if there is $b\in A$ such that
\begin{eqnarray}
 c & = & c(b),\label{eq:surj10}\\
 \pair{\om_i,\nu} & = & \val\om_i(b),\quad\forall i\in \{1,\ldots,n\},\label{eq:surj20}\\
\sg(b) & = & s(b)\,\text{ for some }s\in W,\label{eq:surj30}\text{
and}\\
\val\Om(b) & = & 0\label{eq:surj40}.
\end{eqnarray}
Clearly, \eqref{eq:surj20} holds for $b=a$. For \eqref{eq:surj40}, note that
$$
\al(a)=\eps^{\pair{\al,\nu}}\cdot\al(u)=\al(u),\quad\forall\al\in R,
$$
where on the right we consider $\al$ as a character of $A_{sc}$. Since $u$ is
regular, this gives \eqref{eq:surj40} for $b=a$. Thus, it is enough to show
that $\sg(a)=g(a)$. Since $m_{l+1},\ldots,m_n\in\bZ$ and $W$ acts trivially on
$\om_{l+1},\ldots,\om_n$, this is equivalent to \bbe\label{eq:**}
\sg(\om_i(a))=g^{-1}\om_i(a),\quad\forall i\in\{1,\ldots,l\}.\ee For each $i$
using \eqref{eq:starr} we have
$$
\sg(\om_i(a))=\eps^{m_i}\cdot\exp\left(2\pi\sqrt{-1}m_i\right)\cdot\varpi_i(u)=\eps^{m_i}\cdot
g^{-1}\varpi_i(u)=g^{-1}\om_i(a),
$$
which proves \eqref{eq:**}.

\mbr

Let us show now that Theorem \ref{t:3} holds for $X_{\nu}$. Let $a\in A$ and
$h\in W$ be such that $c(a)\in X_{\nu}$ and $\sg(a)=h(a)$. As was explained in
the last paragraph of the proof of Lemma \ref{l:en}, without loss of generality we
can assume that $\val\la(a)=\pair{\la,\nu}$ for any $\la\in X^*(A)$. We need to
show that $h$ is conjugate to $w(\nu)=\psi_G(\nu)$. Note first that there
exists an element $z$ of the center $Z$ of $G$ such that
\begin{eqnarray*}
&& \om_i(za)=\eps^{m_i}\cdot(v_i+\dotsb),\quad v_i\in\bC^{\times},\,\, i\in\{1,\ldots,n\},\\
&& v_{l+1}=\dotsb=v_n=1.
\end{eqnarray*}
Here the dots in the formula for $\om_i(za)$ denote an element in the maximal ideal of the valuation ring in $\Fb$. Let $u\in T=\Spec\bC[P(R)]$ be given by
$$
\varpi_i(u)=v_i,\quad 1\leq i\leq l.
$$
We claim that $u$ is regular and \bbe\label{eq:***} h(u)=x_{\mu} u,\ee
which together with Proposition \ref{p:sys} and the fact that $\psi_G(\nu)=\psi
(\bar{\mu})$ implies Theorem \ref{t:3}.

\mbr

Since $z\in Z$, we have
$$
\al(a)=\al(za)=\al(u)+\dotsb,\quad \forall\al\in R,
$$
where the dots denote elements in the maximal ideal of the valuation ring in $\Fb$. Since
$\val\Om(a)=0$, we conclude that $\al(u)\ne 1$ for any $\al\in R$, i.e., $u$ is
regular. Furthermore,
for each $i\in\{1,\ldots,l\}$ we get
\begin{eqnarray*}
h^{-1}\om_i(za)&=&\eps^{m_i}\cdot(h^{-1}\varpi_i(u)+\dotsb),\text{ and}\\
\sg(\om_i(za))&=&\eps^{m_i}\cdot\exp(2\pi\sqrt{-1}m_i)\cdot(\varpi_i(u)+\dotsb),
\end{eqnarray*}
where as usual the dots denote elements in the maximal ideal of the valuation ring in $\Fb$. Since
$$
\sg(za)=h(za)\Longleftrightarrow \sg(a)=h(a),
$$ we obtain
$$
h^{-1}\varpi_i(u)=\exp(2\pi\sqrt{-1}m_i)\cdot\varpi_i(u),\quad\forall
i\in\{1,\ldots,l\},
$$
which is equivalent to \eqref{eq:***} by \eqref{eq:starr}.
\end{proof}


\section{Springer's theory for the torus $A$}\label{s:Springer}

This section (which contains analogues of some results of Springer \cite{s})
and the next one can be read independently of the rest of the paper, and are
interesting in their own right. The main results are Theorem \ref{t:abst} and
Theorem \ref{t:s}. The reader who is only interested in Newton stratifications
may decide to take the statement of Theorem \ref{t:s} on faith and skip
directly to \S\ref{ss:Conclusion}, which contains a proof of Proposition
\ref{p:sys}.

\ssect{General conjugacy theorem}\label{ss:abstract}
Let $X$ be a separated algebraic variety over $\bC$, and let $Z$ and $W$ be finite
groups acting on $X$ by morphisms. For $z\in Z$ (resp., $w\in W$) and $x\in X$
we denote by $z\cdot x$ (resp., $w(x)$) the action of $z$ (resp., of $w$) on
$x$. Assume that the actions of $Z$ and $W$ commute, i.e.,
$$
w(z\cdot x)=z\cdot w(x),\quad\forall w\in W,\,z\in Z,\,x\in X.
$$
Denote $X^\circ:=\left\{x\in X\vert\,W_x=\{1\}\right\}$. Let $Y=X/W$ be the
quotient considered as a topological space with the quotient topology of the Zariski topology on $X$. Since
the actions of $Z$ and $W$ on $X$ commute, we have the induced action of $Z$ on
$Y$.

\begin{thm}\label{t:abst}
Let $z\in Z$ be such that $$Y^z:=\left\{y\in Y\vert\,z\cdot y=y\right\}$$ is
irreducible. If there exist $w_1,w_2\in W$ and $x_1,x_2\in X^\circ$ such that
$$
w_i(x_i)=z\cdot x_i,\quad i=1,2,
$$
then $w_1$ is conjugate to $w_2$.
\end{thm}
\begin{proof}
Let $f:X\trar Y$ denote the quotient map. First, note that without loss of
generality we can assume that $X=X^\circ$. Indeed, $X^\circ$ is an open, $W$- and
$Z$-invariant subset of $X$. In particular, $X^\circ$ is separated. Also, by
assumption, $f(X^\circ)^z$ is nonempty and it is an open subset of $Y^z$, because
$f$ is open and $f(X^\circ)^z=f(X^\circ)\cap Y^z$. Thus, $f(X^\circ)^z$ is irreducible.

\sbr

Thus assume that $X=X^\circ$. For $w\in W$ denote $$X(w,z):=\left\{x\in X\vert\,w(x)=z\cdot
x\right\}.$$ Since $X$ is separated, each $X(w,z)$ is closed. Also, for $w\ne
w'\in W$ we have $X(w,z)\cap X(w',z)=\emptyset$. Therefore
$$
f^{-1}(Y^z)=\coprod_{w\in W}X(w,z),
$$
and each $X(w,z)$ is open in $f^{-1}(Y^z)$. Since the restriction of $f$ to
$f^{-1}(Y^z)$ is also an open map, we see that $f\left(X(w_1,z)\right)$ and
$f\left(X(w_2,z)\right)$ intersect, being nonempty open subsets of the irreducible space
$Y^z$. Hence there exist $u,v\in X$ and $g\in W$ such that
$u=g(v)$, $w_1(u)=z\cdot u$, and $w_2(v)=z\cdot v$. Combining these three
equalities and taking into account that $X=X^\circ$, we get $w_1=gw_2g^{-1}$.
\end{proof}

\begin{rem}
This argument was obtained by a careful analysis of the proof of Theorem 4.2(iv) in \cite{s}.
\end{rem}

\begin{cor}\label{c:dim}
Suppose $X$ is a quasi-projective variety $($so that $Y$ has the structure of an
algebraic variety where the underlying topology coincides with the quotient topology$)$, and
$z\in Z$ is such that $Y^z$ is irreducible. If $X^\circ\cap
X(w,z)\ne\emptyset$ for some $w\in W$, then $\dim X(w,z)=\dim Y^z$. Conversely,
assume that $Y^z\cap f(X^\circ)\ne\emptyset$, and let $w\in W$. If $C$ is an irreducible
component of $X(w,z)$ with $\dim C=\dim Y^z$, then $X^\circ\cap C\ne\emptyset$.
\end{cor}

\begin{proof}
Let  $X^\circ\cap X(w,z)\ne\emptyset$. Since $f$ is open and finite, we have
\[
\dim X(w,z) \geq \dim
\left(X^\circ\cap X(w,z)\right)=\dim f(X^\circ\cap X(w,z))=\dim Y^z \geq \dim X(w,z),
\]
where the second equality holds because $f(X^\circ\cap X(w,z))$ is nonempty and open in $Y^z$.

\sbr

Let $\dim C=\dim Y^z$. Then $f(C)$ is dense in $Y^z$, hence intersects $f(X^\circ)$
nontrivially, which implies $X^\circ\cap C\ne\emptyset$, because $X^\circ$ is
$W$-invariant.
\end{proof}

\ssect{New conventions}\label{ss:Conventions} In the rest of the paper we will
be working in a special case of the setup described in \S\ref{ss:setup}. Let
$G$ be a connected and simply connected simple algebraic group over $\bC$; thus
$G_{sc}=G$ and $A_{sc}=A$. All our algebraic varieties (resp., groups) will be
defined over $\bC$, and will be implicitly identified with their sets (resp.,
groups) of complex points; in particular, $G$, $A$, $B$ will stand for
$G(\bC)$, $A(\bC)$, $B(\bC)$, respectively. The tangent space to an algebraic
variety $Y$ at a point $y\in Y$ will be denoted by $T_y Y$, and if $f:Y\rar Y'$
is a morphism to another algebraic variety $Y'$, its differential at $y$ will
be denoted by $D_y f : T_y Y \rar T_{f(y)} Y'$.

\subsection{A variation on Springer's theory of regular elements}\label{ss:springer}
The Weyl group $W=W(G,A)$ acts naturally on $A$ and on the Lie algebra
$\Lie(A)$ of $A$. It is well known that $W$ is a finite complex reflection
group in $\Lie(A)$; thus Springer's results \cite{s} apply to it. In
particular, we recall that Springer introduces the notion of a {\em regular
element} of $W$, namely, it is an element $g\in W$ which has a regular
eigenvector $v\in\Lie(A)$ (the order of the corresponding eigenvalue is then
necessarily equal to the order of $g$). Further, he proves a number of useful
and nontrivial results about regular elements; this includes the fact that if
there exists a regular element $g\in W$ of order $d\in\bN$, then all such
elements form a single conjugacy class ({\em op.~cit.}, Theorem 4.2(iv)), as
well as an explicit determination of the eigenvalues of a regular element ({\em
op.~cit.}, Theorem 4.2(v)).

\mbr

If one tries to find a version of Springer's theory for the action of $W$ on
the torus $A$ itself, rather than on $\Lie(A)$, the first problem one
encounters is how to define the analogue of an eigenvector. Indeed, if $\ze$ is
a root of unity in $\bC$, the equation $g(v)=\ze v$ makes no sense in $A$. Thus
one is led to trying to replace the group of roots of unity by an abelian group
acting on $A$ in a way which commutes with the Weyl group action. We choose the
most naive answer (in some sense), yet one which leads to nontrivial results
that have several interesting applications, as the present paper already
demonstrates. Namely, the group of roots of unity will be replaced by $Z$, the
center of $G$, acting on $A$ by multiplication. We recall that $Z$ is a finite
abelian group which is non-canonically isomorphic to $P(R^\vee)/Q(R^\vee)$, and
we recall the homomorphism $\psi:P(R^\vee)/Q(R^\vee)\rar W$ used earlier in the
paper (its definition is reviewed in \S\ref{ss:psi} below).

\mbr

Our main result in this setup is the following

\begin{thm}\label{t:s}
\begin{enumerate}[$(a)$]
\item For every $x\in Z$, there exists $g\in W$ for which there is
a regular element $u\in A$ satisfying $g(u)=x u$. Moreover, if we choose an
embedding of groups $\bQ/\bZ\hrar\cst$, which induces an isomorphism $Z\cong
P(R^\vee)/Q(R^\vee)$, and denote by $\mu_x\in P(R^\vee)/Q(R^\vee)$ the element
corresponding to $x$ under this isomorphism, then one can take $g=\psi(\mu_x)$.
\item If $x\in Z$ is fixed, the element $g$ satisfying the property above is
determined uniquely up to conjugation.
\item The eigenvalues of such an element $g$ acting on
$\Lie(A)$, counting the multiplicities, are precisely the complex numbers
$\{\vp_i(x)^{-1}\}_{i\in I}$.
\end{enumerate}
\end{thm}

Part (a) of this theorem follows from a stronger result (Theorem
\ref{t:monodromy}) proved in Section \ref{s:geometric}. Namely, it is possible
to choose {\em one} regular element $u\in A$ that works for all $x\in Z$
simultaneously. The rest of the section is devoted to the proofs of parts (b)
and (c).

\subsection{Proof of Theorem \ref{t:s}(b)}\label{ss:Springer} Let $l$ denote the
rank of $G$ and for $x\in Z$ let $a(x)$ denote the number of elements in the
set $\{ 1\leq i\leq l \,\bigl\lvert\, \vp_i(x)=1 \}$. Also, as in
\S\ref{ss:abstract}, we introduce, for each $g\in W$ and $x\in Z$, the subset
\[
A(g,x) = \{ u\in A \,\bigl\lvert\, g(u)=x u \}.
\]
Clearly, Theorem \ref{t:s}(b) follows from the following more
precise proposition, which is an analogue of Theorem 4.2(iv) in \cite{s}:
\begin{prop}
If $g\in W$, $x\in Z$, the following conditions are equivalent:
\begin{enumerate}[$(i)$]
\item $\dim A(g,x)=a(x)$;
\item there exists a regular element $u\in A(g,x)$;
\item there exists a regular element in every connected component of $A(g,x)$.
\end{enumerate}
Moreover, for fixed $x$, the elements of $W$ satisfying these properties form a
single conjugacy class.
\end{prop}

\begin{proof}
Note that although $A(g,x)$ may be disconnected, it is clearly a torsor for the
subgroup $A^g$ consisting of the elements of $A$ that are fixed by $g$. In
particular, $A(g,x)$ is smooth, its connected components ($=$irreducible
components) correspond to the connected components of $A^g$, and consequently
have the same dimension.

\sbr

We recall from \S\ref{ss:adj-quotient} that we have a $W$-invariant polynomial
morphism $c:A\rar\bA=\bA^l$ whose coordinates $c_i$ are given by
\[
c_i(u) = \sum_{\la\in W\vp_i} e^\la(u),\quad u\in A,
\]
and that $c$ identifies $\bA$ with the quotient $A/W$. We apply Theorem
\ref{t:abst} and Corollary \ref{c:dim} to $A$ with the action of $Z$ by
multiplication and the natural action of $W$. We only need to show that $\bA^x$
is irreducible. We have
$$
\bA^x=\left\{(c_i)\in\bA \,\,\bigl\lvert\,\, c_i=\varpi_i(x)\cdot c_i,\,\,
1\leq i\leq l \right\},
$$
which implies that $\bA^x$ is an affine space of dimension $a(x)$. Thus, the
assumptions of Theorem \ref{t:abst} are satisfied. Using Theorem \ref{t:s}(a),
we see that Corollary \ref{c:dim} can be applied to prove the equivalence of (i) -- (iii), and
the proposition follows from these two results.
\end{proof}

\begin{rem}
If $G$ is a finite reflection group in a complex vector space $V$, the quotient
$V/G$ admits a description similar to the description of $A/W$ used above.
Using this description one can prove Theorem 4.2(iv) of \cite{s} in a way which
is completely analogous to the last proof. In the notation of
\S\ref{ss:abstract}, the role of $X$ is played by $V$, the role of $W$ is
played by $G$, and the role of $Z$ is played by the group of roots of unity of
a fixed order $d\in\bN$. In fact, this argument is essentially identical to the
one used by Springer, modulo some elementary simplifications.
\end{rem}

\ssect{Proof of Theorem \ref{t:s}(c)}\label{ss:Eigenvalues} This subsection is
an imitation of the proof of Theorem 4.2(v) in \cite{s}. However, we give a
coordinate-free version of Springer's argument.

\sbr

Note that the map $c:A\rar\bA$ has the property that the differential $D_u
c:T_u A \rar T_{c(u)}\bA$ is an isomorphism if and only if $u\in A$ is regular.
For each $x\in Z$, we define $m_x:A\rar A$ to be the map $u\mapsto xu$; since
$x$ is $W$-invariant, it is clear that we have a commutative diagram
\begin{equation}\label{e:Dilation}
\xymatrix{
  A \ar[dd]_c \ar[rr]^{m_x}   &    &   A \ar[dd]^c \\
    &     &     \\
  \bA  \ar[rr]^{\mb_x} &    &   \bA
   }
\end{equation}
where $\mb_x:\bA\rar\bA$ is defined by $(a_i)_{i\in I}\mapsto (\vp_i(x)\cdot
a_i)_{i\in I}$. In particular, given an element $u\in A$, we have (imitating
\cite{s}, \S3):
\[
xu=g(u)\quad\text{for some  } g\in W \iff c(xu)=c(u) \iff \mb_x(c(u))=c(u)
\]
\[
\iff \text{for each } i\in I, \text{ either } c_i(u)=0, \text{ or } \vp_i(x)=1.
\]

\sbr

Let $g\in G$ and $x\in Z$ be fixed, and let us also fix a regular element $u\in
A(g,x)$ (assuming it exists). The diagram \eqref{e:Dilation} and the fact that
$c\circ g=c$ yield two commutative diagrams (note that $c(xu)=c(g(u))=c(u)$):
\[
\xymatrix{
  T_u A \ar[ddr]_{D_u c} \ar[rr]^{D_u g}   &    &   T_{xu}A \ar[ddl]^{D_{xu}c} \\
    &     &     \\
   &     T_{c(u)}\bA     &
   }
  \qquad\text{and}\qquad
\xymatrix{
  T_u A \ar[dd]_{D_u c} \ar[rr]^{D_u m_x}   &    &   T_{xu}A \ar[dd]^{D_{xu}c} \\
    &     &     \\
  T_{c(u)} \bA  \ar[rr]^{D_{c(u)}\mb_x} &    &  T_{c(u)} \bA
   }
\]
Since $u$, and hence $xu$, are regular, the maps $D_u c$ and $D_{xu} c$ are
isomorphisms. Therefore the automorphism
\[
(D_u g)^{-1}\circ D_u m_x : T_u A \rarab{\simeq} T_u A
\]
is equal to
\begin{eqnarray*}
(D_u g)^{-1}\circ D_u m_x &=& (D_u c)^{-1} \circ (D_{xu} c) \circ (D_{xu} c)^{-1} \circ (D_{c(u)}\mb_x) \circ D_u c \\
 &=& (D_u c)^{-1} \circ (D_{c(u)}\mb_x) \circ (D_u c).
\end{eqnarray*}
In particular, $(D_u g)^{-1}\circ D_u m_x$ has the same eigenvalues as
$D_{c(u)} \mb_x$. But $D_{c(u)} \mb_x$ is obviously given by a diagonal matrix
in the standard basis, whose eigenvalues are the numbers $\vp_i(x)$. On the
other hand, the chain rule yields
\[
(D_u g)^{-1}\circ D_u m_x = D_u (g^{-1}\circ m_x) = D_u(m_x\circ g^{-1}) =
(D_{g^{-1}(u)} m_x)\circ D_u(g^{-1}),
\]
and $(D_{g^{-1}(u)} m_x)\circ D_u(g^{-1})$ has the same eigenvalues as
\[
(D_1 m_u)^{-1} \circ (D_{g^{-1}(u)} m_x) \circ D_u(g^{-1}) \circ D_1 m_u = D_1(
m_u^{-1}\circ m_x \circ g^{-1} \circ m_u) = D_1(g^{-1}).
\]
But $D_1(g^{-1})$, the differential of $g^{-1}$ at the identity element $1\in
A$, is nothing but the automorphism by which $g^{-1}\in W$ acts on $\Lie(A)$.
Hence we have shown that the eigenvalues of $g^{-1}$ are the numbers
$\vp_i(x)$, completing the proof of part (c) of Theorem \ref{t:s}.

\ssect{Proof of Proposition \ref{p:sys}}\label{ss:Conclusion} It is clear that
if Proposition \ref{p:sys} holds for two root systems, then it also holds for
their direct sum; thus it suffices to prove it for a reduced and irreducible
root system. It is then obvious that parts (a) and (b) of Theorem \ref{t:s} imply
the desired result, which finally completes the proof of Theorem \ref{t:3}, and hence of Theorem \ref{t:1}.

\begin{rem}\label{r:kottwitz}
As an additional bonus, we note that Theorem \ref{t:s} yields an alternate
proof of Lemma 4.1.1 of \cite{k}, namely, a geometric one which does not use
the classification of simple Lie algebras over $\bC$. Indeed, as explained on
p.~14 of {\em op.~cit.}, the lemma reduces to the statement that, with our
notation, if we choose the embedding $\iota:\bQ/\bZ\hrar\cst$ defined by
$\iota(q)=\exp(-2\pi\sqrt{-1}\cdot q)$ (this is the same choice as the one that
was made in Section \ref{s:systems}, but is complex conjugate to the one used
by Kottwitz in \cite{k}), and the corresponding isomorphism $Z\cong
P(R^\vee)/Q(R^\vee)$, then the representation of $Z$ on $\Lie(A)^*$ obtained
from the natural representation of $W$ via the homomorphism $Z\rar W$,
$x\mapsto\psi(\mu_x)$, decomposes into the direct sum of the $1$-dimensional
representations given by the characters $\vp_i\bigl\lvert_Z$. However, basic
character theory of finite groups implies that it is enough to show that for
each $x\in Z$, the eigenvalues of $\psi(\mu_x)$ on $\Lie(A)^*$ are the numbers
$\vp_i(x)$, counting the multiplicities, and this is precisely the statement
dual to part (c) of Theorem \ref{t:s} above, in view of part (a) of Theorem
\ref{t:s}.
\end{rem}


\section{Geometric interpretation of $\psi$}\label{s:geometric}

We keep the conventions of \S\ref{ss:Conventions} and use the notation ($G$,
$A$, $B$, $R$, $W$, etc.) introduced in \S\ref{ss:setup}. In this section we
study in more detail the homomorphism $\psi:P(R^\vee)/Q(R^\vee)\to W$
(determined by the choice of $B$) and prove part (a) of Theorem \ref{t:s}. The
main result of the section is Theorem \ref{t:monodromy}, which provides a
geometric definition of $\psi$.

\ssect{Classical definition of the homomorphism $\psi$}\label{ss:psi} We begin
by recalling the definition of the homomorphism $\psi$ that appears in
\cite{b}. Perhaps we provide a little more information than is strictly needed
for our purposes, but our presentation is somewhat clearer than Bourbaki's, and
the additional facts that we mention may help the reader understand the general
picture.

\mbr

Let $V=Q(R)\tens_\bZ\bR$; note that our notation is different from the one used
in \cite{k}, since Kottwitz denotes by $V$ the complex space
$Q(R)\tens_\bZ\bC$. With our notation, $R$ can be thought of as an abstract
reduced and irreducible root system in the real vector space $V$ (\cite{b},
Chapter VI, \S1, no.~1). We introduce the corresponding affine root system
$\Rt$ in an ad hoc manner as follows. The root lattice of $\Rt$ is by
definition the lattice $Q(\Rt):=Q(R)\oplus\bZ\cdot\al_0$ in the vector space
$\Vt:=V\oplus\bR\cdot\al_0$, where $\al_0$ is just an auxiliary symbol. We call
$\al_0$ the {\em affine simple root}. Let $\alt$ denote the highest root of
$R$, and define $\de=\al_0+\alt\in Q(\Rt)$. We set
\[
\Rt = \bigl\{ n\de + \al \,\lvert\, \al\in R,\ n\in\bZ \bigr\} \bigcup \bigl\{
m\de \,\lvert\, m\in\bZ\setminus\{0\}\}.
\]
This is the {\em affine root system} associated to $R$. The elements of $\Rt$
of the form $m\de$, $m\in\bZ\setminus\{0\}$, are called {\em imaginary roots};
all other elements are called {\em real} ({\em affine}) {\em roots}. A real
root $n\de+\al$ is called {\em Dynkin} if $n=0$, and {\em non-Dynkin}
otherwise. Thus $R$ is identified with the subset of $\Rt$ consisting of Dynkin
roots.

\mbr

The elements $\{\al_i\}_{i\in I\cup \{0\}}$ are the {\em simple roots} of the
system $\Rt$. Thus $\al_0$ is the unique non-Dynkin simple root. The set
$I\cup\{0\}$ is in natural bijection with the set of vertices of the extended
Dynkin diagram corresponding to the root system $R$. Moreover, let us write
\[
\alt = \sum_{i\in I} \de_i \al_i, \qquad \de_i\in\bN,
\]
(note that Bourbaki uses a different notation: $\alt=\sum n_i\al_i$), and set
$\de_0=1$. Then we have (by construction)
\[
\de = \sum_{i\in I\cup\{0\}} \de_i \al_i.
\]

\mbr

We identify $V$ with a quotient of $\Vt$, namely, $\Vt\bigl/(\bR\cdot\de)\cong
V$.  Explicitly, this isomorphism is the inverse of the composition of the
natural inclusion $V=V\oplus(0)\hrar\Vt$, followed by the projection. This
isomorphism takes the image of $\al_i\in\Vt$ in $\Vt\bigl/(\bR\cdot\de)$ to
$\al_i\in V$ for all $i\in I$, and it takes the image of $\al_0\in\Vt$ in
$\Vt\bigl/(\bR\cdot\de)$ to $-\alt\in V$. Dually, we identify $V^*$ with a
subspace of $\Vt^*$:
\[
V^* \cong \{ f\in \Vt^* \,\bigl\lvert\, f(\de)=0 \}.
\]
In particular,
\begin{equation}\label{e:Rvee}
P(R^\vee) \cong \{ f\in\Vt^* \,\bigl\lvert\, f(R)\subseteq\bZ,\ \ f(\de)=0 \}.
\end{equation}
We also define
\[
E = \{ f\in \Vt^* \,\bigl\lvert\, f(\de)=1 \}.
\]
This is clearly an affine space for the vector space $V^*$. We now define a
linear action $\mu\mapsto t_\mu$ of $P(R^\vee)$ on $\Vt^*$ as follows: if
$\mu\in P(R^\vee)$ and $f\in\Vt^*$, then $t_\mu(f)=f+f(\de)\cdot\mu$. Note that
this action preserves $E$, thanks to \eqref{e:Rvee}, and the induced action of
$P(R^\vee)$ on $E$ is simply by translations. We also define $\mu\mapsto
t^*_\mu$ to be the contragredient action of $P(R^\vee)$ on $\Vt$. On the other
hand, we define an action of the Weyl group $W$ on $\Vt$ by declaring that $W$
acts trivially on $\de$ and acts by its usual action on $V$. Then we also
obtain the contragredient action of $W$ on $\Vt^*$, and it is easy to check
that the action of $W$ on $\Vt$ (resp., on $\Vt^*$) normalizes the action of
$P(R^\vee)$ on $\Vt$ (resp., on $\Vt^*$); in fact, e.g., if $w\in W$ and
$\mu\in P(R^\vee)$, then $wt^*_\mu w^{-1} = t^*_{w(\mu)}$. This equation
implies that we have an action of the {\em extended Weyl group}
$W_e:=P(R^\vee)\rtimes W$ on $\Vt$ and $\Vt^*$, and, moreover, this action
preserves $E$. Note also that $W_e$ contains as a subgroup the {\em affine Weyl
group} $W_a:=Q(R^\vee)\rtimes W$.

\mbr

The last key ingredient is the following fact. For every real affine root
$\al\in\Rt$, we have the corresponding root hyperplane
$\Ker(\al)\subseteq\Vt^*$, and the intersection $\Ker(\al)\cap E$ is a
hyperplane in $E$. The complement of the union of all such affine root
hyperplanes in $E$ is a disjoint union of bounded connected open subsets of
$E$, called the {\em alcoves}. Moreover, it is clear that the action of $W_e$
permutes the real affine roots, hence also permutes the alcoves. A fundamental
result (\cite{b}, Chapter VI, \S2, no.~1) is that the affine Weyl group $W_a$
already acts {\em simply transitively} on the alcoves. This implies that if $C$
is a fixed alcove and $\Ga_C$ is the group of automorphisms of $C$ in $W_e$,
then $\Ga_C$ projects isomorphically onto the quotient $W_e/W_a\cong
P(R^\vee)/Q(R^\vee)$. Furthermore, we do have a canonical choice for $C$: it is
the so-called {\em fundamental alcove}, defined by
\[
C = \bigl\{ v\in E \,\,\bigl\lvert\,\, 0<v(\al_i)<1 \ \forall\,i\in I\cup\{0\}
\bigr\}.
\]
If $\xi:\Ga_C\rar P(R^\vee)/Q(R^\vee)$ is the corresponding isomorphism and
$\pr:W_e\rar W$ is the natural projection, then the sought-after homomorphism
$\psi:P(R^\vee)/Q(R^\vee) \rar W$ is defined as the composition
$\psi=\pr\circ\xi^{-1}$. Note that the construction of $\psi$ is completely
canonical once the basis of simple roots $\{\al_i\}_{i\in I}$ has been chosen.

\mbr

We now state the main property of the homomorphism $\psi$ that will be used
below. Let us denote by $\{\ob_i\}_{i\in I}$ the basis of $P(R^\vee)$ dual to
the basis $\{\al_i\}_{i\in I}$ of $Q(R)$. For notational convenience, we set
$\ob_0=0$. Following Bourbaki (\cite{b}, Chapter VI, \S2, no.~2), we define a
subset $J\subseteq I$ by
\[
J = \{ i\in I \,\bigl\lvert\, \de_i=1 \}.
\]
For example, $J=I$ if $R$ is of type $A$, and $J\subsetneq I$ otherwise.
\begin{lem}\label{l:psi} The elements $\{\ob_i\}_{i\in J\cup\{0\}}$ form a
complete set of representatives for the quotient $P(R^\vee)/Q(R^\vee)$.
Moreover, for each $i\in J\cup\{0\}$, we have
$\xi^{-1}\left(\left[\ob_i\right]\right)=t_{\ob_i}\circ\psi(\ob_i)$, where
$\left[\ob_i\right]$ denotes the image of $\ob_i$ in $P(R^\vee)/Q(R^\vee)$, and
$\psi(\ob_i)$ is viewed as an element of $W_e$ via the obvious inclusion
$W\hrar W_e$.
\end{lem}
\begin{proof}
See \cite{b}, Chapter VI, \S2, no.~3, Proposition 6 and
its Corollary.
\end{proof}

\ssect{Interpretation of $\psi$ in terms of monodromy}\label{ss:monodromy} The
homomorphism $\psi$ described above does not appear to be very well understood,
since the purely algebraic definition of \cite{b} is not very enlightening. The
results of this section and the previous one are a first step towards
understanding the homomorphism $\psi$ in more concrete terms. One drawback of
Theorem \ref{t:s} is that it only describes the conjugacy classes of each of
the elements $\psi(\mu)$ separately, rather than the conjugacy class of the
whole homomorphism $\psi$. This drawback is corrected in Theorem
\ref{t:monodromy} below.

\mbr

We now concentrate on the question of whether the homomorphism $\psi$ admits a
geometric, rather than algebraic, description. One way to state this question
precisely is as follows: do there exist a smooth complex algebraic variety $Y$
whose fundamental group (with respect to the complex topology) is isomorphic to
$P(R^\vee)/Q(R^\vee)$ and a $W$-torsor $X$ over $Y$ such that the induced
monodromy homomorphism $P(R^\vee)/Q(R^\vee)\rar W$ is conjugate to $\psi$? We
do not know a definite answer to this very question, but we do have a
``semi-positive'' answer, which may be sufficient for many practical purposes.
The formulation of the answer is due to Ilya Zakharevich. The main idea is as
follows.

\mbr

Let $G_{ad}=G/Z$ be the simple adjoint group corresponding to $G$, let
$G^{rs}\subset G$ and $G^{rs}_{ad}\subset G_{ad}$ be the (Zariski open) subsets
consisting of the regular semisimple elements of the groups $G$ and $G_{ad}$,
and consider the following diagram:
\begin{equation}\label{e:coverings}
\xymatrix{
  \ar[dd]_{p'} \widetilde{G_{ad}^{rs}}  &     & \ar[ll]_{q'} \ar[dd]^p \widetilde{G^{rs}}   \\
     &        &     \\
   G_{ad}^{rs}   &     &   \ar[ll]_q G^{rs}
   }
\end{equation}
Here $\widetilde{G^{rs}}$ denotes the incidence variety of pairs $(g,B')$
consisting of an element $g\in G^{rs}$ and a Borel subgroup $B'\subset G$
containing $g$, and $p$ is the obvious projection. It is well known that $p$ is
a $W$-torsor, and, in particular, topologically it is a covering map. The
$W$-torsor $p':\widetilde{G^{rs}_{ad}}\rar G^{rs}_{ad}$ is defined in a similar
manner; $q$ is induced by the quotient map $G\rar G_{ad}$, and $q'$ is induced
by $q$ in the obvious way. All maps in the diagram above are surjective, and
$q$ can also be thought of as a $Z$-torsor.

\mbr

The monodromy of $p'$ yields a homomorphism $\tau:\pi_1(G_{ad}^{rs})\rar W$,
which is well defined up to conjugation. On the other hand, $q$ induces an
exact sequence of groups
\begin{equation}\label{e:fund-grps}
\pi_1(G^{rs}) \rarab{q_*} \pi_1(G^{rs}_{ad}) \rar Z \rar 1.
\end{equation}
Of course, if the composition $\tau\circ q_*$ were trivial, then $\tau$ would
induce a homomorphism $Z\rar W$. However, this is far from being true, for the
triviality of $\tau\circ q_*$ would imply that $p$ is a trivial $W$-torsor,
which is certainly false. Instead, we proceed as follows. Let us choose a
compact form $K\subset G$ of $G$, and write $K_{ad}=K/Z$, which is a compact
form of $G_{ad}$. The diagram \eqref{e:coverings} can be restricted to the
``compact part'', which produces
\begin{equation}\label{e:coverings-compact}
\xymatrix{
  \ar[dd]_{p'_K} \widetilde{K_{ad}^{rs}}  &     & \ar[ll]_{q'_K} \ar[dd]^{p_K} \widetilde{K^{rs}}   \\
     &        &     \\
   K_{ad}^{rs}   &     &   \ar[ll]_{q_K} K^{rs}
   }
\end{equation}
Of course, $K^{rs}$ (resp., $K^{rs}_{ad}$) just stands for the set of regular
elements of $K$ (resp., $K_{ad}$), since all elements of $K$ and $K_{ad}$ are
automatically semisimple. The advantage of considering this diagram lies in the following
\begin{lem}\label{l:codim}
The complement $K\setminus K^{rs}$ has real codimension $\geq 3$ in $K$.
\end{lem}
\begin{proof}
Let $T=A\cap K$, which is a maximal torus of $K$. Let $T^{reg}=K^{rs}\cap T$
denote the set of regular elements of $T$. It is well known that every element
$y\in K$ is conjugate to an element of $T$, and $y$ is regular if and only if
it is conjugate to an element of $T^{reg}$. Thus the conjugation map
\[
K \times \bigl(T\setminus T^{reg}\bigr) \rar K\setminus K^{rs}
\]
is surjective. Moreover, its fibers have dimension $\geq T+2$, because if $t\in
T$ is not regular, its centralizer in $K$ contains $T$ as well as a nontrivial
compact semisimple Lie group. Hence
\[
\dim(K\setminus K^{rs}) \leq \dim(K) + \dim(T\setminus T^{reg}) - \dim(T) -2 \leq \dim(K)-3,
\]
completing the proof.
\end{proof}

A standard fact in algebraic topology now yields

\begin{cor}\label{c:fund}
The inclusions $K^{rs}\hrar K$ and $K_{ad}^{rs}\hrar K_{ad}$ induce isomorphisms on the fundamental groups.
\end{cor}

Thus $K^{rs}$ is simply connected, and the natural map $\pi_1(K^{rs}_{ad})\rar
Z$ is an isomorphism. In particular, the monodromy of $p'_K$ does determine a
conjugacy class of homomorphisms $\tau_K : Z \rar W$. Our last main result is

\begin{thm}\label{t:monodromy}
In the situation of Theorem $\ref{t:s}$, there exists a regular element $u\in
A\cap K$ such that $\psi(\mu_x)(u)=x u$ for all $x\in Z$. If $\overline{u}$
denotes its image in $K_{ad}^{rs}$, then the corresponding monodromy
homomorphism $\tau_K:Z\cong\pi_1(K_{ad}^{rs},\overline{u})\rar W$ is conjugate
to the homomorphism $x\mapsto\psi(\mu_x)$.
\end{thm}

Note that this theorem is not quite what we were looking for, since the
construction is based on the compact form $K$, and hence cannot be formulated
in a completely algebro-geometric way. The situation can be improved slightly
by noting that the inclusion $K^{rs}_{ad}\hrar G^{rs}_{ad}$ induces a
homomorphism
\[
j : Z\rarab{\simeq}\pi_1(K^{rs}_{ad})\rar\pi_1(G^{rs}_{ad})
\]
which is a splitting of \eqref{e:fund-grps}. Hence the homomorphism
$x\mapsto\psi(\mu_x)$ is also conjugate to the composition of the monodromy of
the algebraic covering $p':\widetilde{G^{rs}_{ad}}\rar G^{rs}_{ad}$ with $j$.
Unfortunately, the definition of $j$ still involves $K$; however, if we could
characterize abstractly all splittings $j$ of \eqref{e:fund-grps} that arise in
this way, then we would have a manageable algebro-geometric description of the
homomorphism $\psi:P(R^\vee)/Q(R^\vee)\rar W$.

\subsection{Proof of Theorem \ref{t:monodromy}}\label{ss:proof-monodromy} The proof
of the first statement of Theorem \ref{t:monodromy} that we give below was
explained to us by R.~Kottwitz (it is much faster and more transparent than our
original argument). A similar idea was suggested independently by J.-P.~Serre.

\mbr

To simplify notation, we assume that the embedding $\iota:\bQ/\bZ\hrar\cst$ is
chosen to be the same as the one used in Section \ref{s:systems}, namely,
$\iota(q)=\exp(-2\pi\sqrt{-1}\cdot q)$. Let $T=A\cap K$, as before. There is a
natural identification $T=\Hom(P(R),S^1)$. If $V$, $\widetilde{V}$ and
$E\subset \widetilde{V}^*$ are as in \S\ref{ss:psi}, we define an ``exponential
map''
\[
\Exp : E \rar T
\]
by the formula
\[
\Exp(f)(v)=\exp\bigl(2\pi\sqrt{-1}\cdot f(v)\bigr) \qquad \forall\, v\in P(R);
\]
it identifies $T$ with the quotient of $E$ by the translation action of
$Q(R^\vee)$. In particular, let $e\in C$ denote the barycenter of the
fundamental alcove $C\subset E$ introduced in \S\ref{ss:psi} (i.e., $e$ is the
arithmetic mean of the vertices of $C$). Since $e$ does not lie in any of the
affine root hyperplanes in $E$, the element $u:=\Exp(e)\in T$ is regular. On
the other hand, $\Ga_C$ fixes $e$ because it preserves $C$. According to Lemma
\ref{l:psi}, this means that $t_{-\nu}(e)\equiv\psi(\nu)(e)$ modulo $Q(R^\vee)$
for all $\nu\in P(R^\vee)/Q(R^\vee)$, which is equivalent to the property
$xu=\psi(\mu_x)(u)$ for all $x\in Z$. This proves the first statement of
Theorem \ref{t:monodromy}.

\mbr

We now prove the second statement. Let $u$ be as above, and let $\overline{u}$
denote the image of $u$ in $K^{rs}_{ad}$. For each $x\in Z$, choose a
continuous path $\ga_x$ from $u$ to $xu$ in $K^{rs}$, and let
$\overline{\ga}_x$ denote its image in $K^{rs}_{ad}$, which is a loop at
$\overline{u}$. Then the map $Z\rar \pi_1(K^{ad}_{rs},\overline{u})$ given by
$x\mapsto[\overline{\ga}_x]$ is a group isomorphism by Corollary \ref{c:fund}.
Since $K^{rs}$ is simply connected, the projection $p_K:\widetilde{K}^{rs}\to
K^{rs}$ has a (unique) continuous section $\sg_K:K^{rs}\to\widetilde{K}^{rs}$
satisfying $\sg_K(u)=(u,B)$, where $B$ is our chosen Borel subgroup of $G$
containing $A\ni u$. The second statement of Theorem \ref{t:monodromy} is
equivalent to the assertion that $\sg_K(xu)=(xu,B^{\psi(\mu_x)})$ for all $x\in
Z$, where $B^{\psi(\mu_x)}$ is the Borel subgroup obtained by conjugating $B$
with the Weyl group element $\psi(\mu_x)\in W$. However, $xu=\psi(\mu_x)(u)$ by
construction, so the equality $\sg_K(xu)=(xu,B^{\psi(\mu_x)})$ follows from the
stronger assertion that $\sg_K(gug^{-1})=(gug^{-1},B^g)$ for all $g\in K$. To
verify the last equality, observe that it holds for $g=1$ by assumption, and
therefore holds for all $g\in K$ by connectedness. This proves Theorem
\ref{t:monodromy}.

\end{document}